\documentclass[final,leqno,onefignum]{siamltex1213}
\usepackage{amssymb}
\usepackage{amsmath}
\usepackage{algorithm}
\usepackage{algpseudocode}
\usepackage{multirow}
\usepackage{comment}
\usepackage{graphicx}
\usepackage{array}
\newcolumntype{C}[1]{>{\centering\let\newline\\\arraybackslash\hspace{0pt}}m{#1}}

\title{A Preconditioned Low-Rank Projection Method with a Rank-Reduction Scheme for Stochastic Partial Differential Equations\thanks{This work was supported by the U.S. Department of Energy Office of Advanced Scientific Computing Research, Applied Mathematics program under award DE-SC0009301 and by the U.S. National Science Foundation under grant DMS1418754.}} 

\author{Kookjin Lee\thanks{Department of Computer Science, University of Maryland, College Park, MD 20742 (\email{klee@cs.umd.edu}).}\and Howard C. Elman\thanks{Department of Computer Science and Institute for Advanced Computer Studies, University of Maryland, College Park, MD 20742 (\email{elman@cs.umd.edu})}}

\begin{document}
\maketitle
\slugger{mms}{xxxx}{xx}{x}{x--x}

\begin{abstract}
In this study, we consider the numerical solution of large systems of linear equations obtained from the stochastic Galerkin formulation of stochastic partial differential equations. We propose an iterative algorithm that exploits the Kronecker product structure of the linear systems. The proposed algorithm efficiently approximates the solutions in low-rank tensor format. Using standard Krylov subspace methods for the data in tensor format is computationally prohibitive due to the rapid growth of tensor ranks during the iterations. To keep tensor ranks low over the entire iteration process, we devise a rank-reduction scheme that can be combined with the iterative algorithm. The proposed rank-reduction scheme identifies an important subspace in the stochastic domain and compresses tensors of high rank on-the-fly during the iterations. The proposed reduction scheme is a multilevel method in that the important subspace can be identified inexpensively in a coarse spatial grid setting. The efficiency of the proposed method is illustrated by numerical experiments on benchmark problems.

\end{abstract}



\pagestyle{myheadings}
\thispagestyle{plain}
\markboth{K. LEE AND H.C. ELMAN}{LOW-RANK PROJECTION METHOD WITH RANK-REDUCTION}

\section{Introduction}
Consider the stochastic elliptic boundary value problem, to find a random function, $u(\pmb{x},\,\xi):  \bar{D} \times \Gamma \rightarrow \mathbb{R}$ that satisfies
\begin{equation}\label{model_problem_param}
\mathcal{L}(a(\pmb{x},\,\xi))(u(\pmb{x},\,\xi)) = f(\pmb{x}) \quad \text{in } D \times \Gamma, \end{equation}
where $\mathcal{L}$ is a linear elliptic operator and $a(\pmb{x},\xi)$ is a positive random field parameterized by a set of random variables $\xi = \{\xi_1,\,\dots,\,\xi_M\}$. The problem is posed on a bounded domain $D \subset \mathbb{R}^2$ with appropriate boundary conditions. Such problems arise, for example, from fluid flow and the transport of chemicals in flows in heterogeneous porous media, where the permeability coefficient is modeled as a random field \cite{ghanem2003stochastic, sun2002stochastic}. 

As the solution method for \eqref{model_problem_param}, we consider the stochastic Galerkin method   \cite{babuvska2002solving, babuska2004galerkin, ghanem2003stochastic}, which, after 
 discretization, leads to a large coupled deterministic system
\begin{equation}
\label{linear_sys}
Au = f,
\end{equation}
for which computations will be expensive for large-scale applications. When the coefficient $a(\pmb{x},\xi)$ has an affine structure depending on a finite number of random variables, the system matrix $A$ can be represented by a sum of Kronecker products of smaller matrices. Matrix operations such as matrix-vector products that take advantage of the tensor format can be performed efficiently, which makes the use of iterative solvers attractive. In this study, we develop a new efficient iterative solver for systems represented in the Kronecker-product structure.

In recent years, many authors started to explore the Kronecker-product structure of such problems and developed iterative algorithms that exploit the structure to reduce computational efforts \cite{ballani2013projection, kressner2010krylov, kressner2011low, matthies2012solving, powell2015efficient, schwab2011sparse}. In addition, it has been shown that the solution of \eqref{model_problem_param} in the stochastic Galerkin setting can be approximated by a tensor of low rank, which further reduces computational effort \cite{benner2013low}. If Krylov subspace methods are used to compute such a solution, however, it may happen that approximate solutions or other auxiliary terms obtained during the course of an iteration do not have low rank, and rank-reduction schemes are required to keep costs under control.

In this study, we will explore a variant of the generalized minimum residual (GMRES) method combined with a rank-reduction strategy that exploits specific features of the stochastic Galerkin formulation. The  strategy we propose is a multilevel scheme that first identifies a low-dimensional subspace, obtained from a coarse-grid spatial discretization, on which a low-rank coarse-grid tensor solution is computed. This solution can be used to estimate the rank of the tensor solution for the desired fine-grid solution. This information is used to define a strategy for rank reduction to be used with iteration on the fine grid space. We show that this strategy enhances the efficiency of preconditioned GMRES for computing the solution.


An outline of the paper is as follows. In section \ref{sec_sg_model}, we review the stochastic Galerkin method and present the Kronecker-product structure of Galerkin systems. 
In section \ref{sec_lr_proj}, we present a preconditioned  projection method for computing approximate solutions in low-rank tensor format. In section \ref{sec_trunc}, we review the conventional approaches and propose a multilevel rank-reduction scheme, which is the main contribution of this work. In section \ref{sec_num_exp}, we illustrate the effectiveness of the low-rank projection method combined with the proposed truncation scheme by numerical experiments on benchmark problems. Finally, in section \ref{sec_con}, we draw some conclusions.

\section{Model problems with random inputs} \label{sec_sg_model}
Consider the steady-state stochastic diffusion equation with homogeneous Dirichlet boundary conditions, 
\begin{equation}
\left\{
\begin{array}{r l l}
-\nabla \cdot (a(\pmb{x},\,\omega) \nabla u(\pmb{x},\,\omega) ) &= f(\pmb{x}) &\text{ in } D \times \Omega,\vspace{2mm}\\
u(\pmb{x},\,\omega) &= 0 &\text{ on } \partial D \times \Omega,
\end{array}
\right.
\label{strong_diff}
\end{equation}
where the diffusion coefficient $a(\pmb{x},\,\omega)$ is a random field and $\omega$ is an elementary event in a probability space. The gradient operator $\nabla$ only acts on the physical domain $D$. The weak form of \eqref{strong_diff} is  to find  $u$ in  $V = H_0^1(D) \otimes L^2(\Omega)$ such that
\begin{equation}
\left\langle \int_D a(\pmb{x},\,\omega) \nabla u(\pmb{x},\,\omega) \cdot \nabla v(\pmb{x},\,\omega) d\pmb{x}    \right\rangle = \left \langle  \int_D f v(\pmb{x},\,\omega) \right \rangle, \quad \forall v(\pmb{x},\,\omega)  \in V  \label{cont_weak_diff}
\end{equation}
where $\langle \cdot \rangle$ refers to expected value with respect to the probability measure on $L_2(\Omega)$ and $V$ is equipped with the gradient norm 
\begin{equation}
\| v \|_V = \int_\Omega \int_D a(\pmb{x},\,\xi) | \nabla v(\pmb{x},\,\xi) |^2\: d\pmb{x}\, dP(\omega).
\end{equation}
If $a(\pmb{x},\,\omega)$ is bounded and strictly positive, 
\begin{equation}
0 < a_{\min} \leq a(\pmb{x},\,\omega) \leq a_{\max} < +\infty, \quad \text{a.e. in } D \times \Omega,
\label{coer_assumption}
\end{equation}
then the Lax-Milgram lemma can be applied to establish existence and uniqueness of a solution $u(\pmb{x},\,\omega) \in V$ of the  variational problem \eqref{cont_weak_diff}.  For the random field $a(\pmb{x},\,\omega)$ with mean $a_0$ and variance $\sigma^2$, we consider a truncated Karhunen-Lo\'eve expansion \cite{loeve1978probability},
\begin{equation}
a(\pmb{x},\, \omega)  = a_0 + \sigma \sum_{i=1}^M \sqrt{\lambda_i} a_i(\pmb{x}) \xi_i(\omega),\label{kl_exp}
\end{equation}
where  $(\lambda_i,\,a_i)$ is an eigenpair of the covariance kernel $C(\pmb{x},\,\pmb{y})$ of the random field. In the following, we use the notation $a(\pmb{x},\,\xi)$  for the random field, let $\Gamma$ denote the image space $\xi(\omega)$, and refer to the expected value of $v(\xi)$ as $\langle v(\xi) \rangle_\rho = \int_\Gamma v(\xi) \rho(\xi) d\xi$.

\subsection{Stochastic Galerkin method}
The stochastic Galerkin method \cite{babuvska2002solving, babuska2004galerkin,ghanem2003stochastic} seeks a finite-dimensional solution $u_{hp}(\pmb{x},\,{\xi}) \in W^h = X_h \otimes S_M$ such that
\begin{equation}
\left\langle \int_D a(\pmb{x},\,{\xi}) \nabla u_{hp}(\pmb{x},\,{\xi}) \cdot \nabla v(\pmb{x},\,{\xi}) d\pmb{x}    \right\rangle_\rho = \left \langle  \int_D f v(\pmb{x},\,{\xi}) \right \rangle_\rho, \quad v(\pmb{x},\,{\xi})  \in W^h  \label{galerkin_form}
\end{equation}
where $X_h$ and $S_M$ are finite-dimensional subspaces of $H_0^1(D)$ and $L_\rho^2(\Gamma)$,
\begin{equation*}
\label{fem_sg}
X_h = \text{span}\{ \phi_r(\pmb{x}) \}_{r=1}^{n_x} \subset H_0^1(D), \quad S_M = \text{span}\{\psi_s ({\xi}) \}_{s=1}^{n_\xi} \subset L_\rho^2(\Gamma),
\end{equation*}
and
\begin{equation}
u_{hp}(\pmb{x},\,{\xi}) = \sum_{s=1}^{n_\xi} \sum_{r=1}^{n_x} u_{r,\,s} \phi_r(\pmb{x}) \psi_s({\xi}). \label{sg_sol}
\end{equation}
Here, $\{ \phi_r \}$ are standard finite element basis functions and $\{ \psi_s \}$ are basis functions for the generalized polynomial chaos expansion(PCE) \cite{xiu2002wiener} consisting of products of orthonormal univariate polynomials: 
$\psi_s({\xi}) = \psi_{{\alpha}({s})}({\xi}) = \prod_{i=1}^M \pi_{\alpha_i(s)} (\xi_i)$
where $\{\pi_{\alpha_i(s)} (\xi_i)\}_{i=1}^M $ are univariate polynomials, and ${\alpha}(s) = (\alpha_1(s),\,\dots\,\alpha_M(s)) \in \mathbb{N}^{M}_0$ is a multi-index where $\alpha_i$ represents the degree of a polynomial in $\xi_i$. 

Let us define an index set $\Lambda_M = \{ {\alpha}(s) \in \mathbb{N}_0^M\}  \label{indexset}$ where each element in $\Lambda_M$ has one-to-one mapping with a natural number. In this study, we set $\Lambda_M$ to be the Total Degree (TD) space $\Lambda_{M,\, p}$:  
\begin{equation}
\Lambda_{M,\, p} = \{{\alpha}(s) \in \mathbb{N}^M_0: \|{\alpha}(s)\|_0 \leq M,\, \|{\alpha}(s)\|_1 \leq p\},  \label{td_set}
\end{equation}
where $\| {\alpha}(s) \|_0$ is the cardinality of a set $\{ k: \alpha_k(s) \neq 0\}$, $\|{\alpha}(s)\|_1 = \sum_{k=1}^M \alpha_k(s)$, and $p$ defines the maximal degree of $\{\psi_i\}_{i=1}^{n_\xi}$. Then the number of PCE basis functions is  $n_\xi = \dim (\Lambda_{M,\,p}) = \frac{(M + p)!}{M!p!}$.

\subsection{Stochastic Galerkin formulation in tensor notation} 
It follows from \eqref{galerkin_form} and \eqref{sg_sol} that the linear system $Au=f$ of \eqref{linear_sys} can be represented in tensor product notation \cite{powell2009block},
\begin{align}
\label{linear_sys_tensor}
\left(G_0 \otimes K_0 + \sum_{l=1}^M G_l \otimes K_l \right) u = g_0 \otimes f_0
\end{align}
where $\otimes$ is the Kronecker product, $\{K_i\}$ are weighted stiffness matrices defined via
\begin{align*}
[K_0]_{ij} = \int_D a_0 \nabla \phi_i (\pmb{x}) \nabla \phi_j (\pmb{x}) d\pmb{x}, \quad
[K_l]_{ij} = \int_D  \tilde{a}_l(\pmb{x}) \nabla \phi_i (\pmb{x}) \nabla \phi_j (\pmb{x}) d\pmb{x}, \:\: l = 1,\,\dots,\,M,
\end{align*} 
$\tilde{a}_l(\pmb{x}) = \sigma \sqrt{\lambda_l}\, a_l(\pmb{x})$, $\{G_i\}$ are ``stochastic'' matrices defined via
\begin{align}
\label{stoch_matrices}
[G_0]_{ij} = \left\langle   \psi_i ({\xi})  \psi_j ({\xi}) \right\rangle_\rho, \quad [G_l]_{ij} =  \left\langle  \xi_l\, \psi_i ({\xi})  \psi_j ({\xi}) \right\rangle_\rho, \quad l = 1,\,\dots,\,M,
\end{align} 
and the vectors $f_0$ and $g_0$ are defined via \begin{align*}
[f_0]_i = \int_D f \phi_i(\pmb{x}) d\pmb{x},\quad [g_0]_i = \left\langle  \psi_i ({\xi}) \right\rangle_\rho.
\end{align*}
Here the Kronecker product of two matrices $G \in \mathbb{R}^{n_\xi \times n_\xi}$ and $K \in \mathbb{R}^{n_x \times n_x}$ is 
\begin{equation*}
G \otimes K = 
\begin{bmatrix}{} 
[G]_{11} K & \dots  & [G]_{1n_\xi}K\\ 
\vdots & & \vdots\\
[G]_{n_\xi1} K & \dots & [G]_{n_\xi n_\xi}K
\end{bmatrix}.
\end{equation*}
Note that $\{ G_l \} _{l=0}^M$ of \eqref{stoch_matrices} are highly sparse because of the orthogonality properties of the stochastic basis functions \cite{ernst2010stochastic}. 

We will make use of an isomorphism between $\mathbb{R}^{n_x  n_\xi}$ and $\mathbb{R}^{n_x \times n_\xi}$ determined by the operators vec($\cdot$) and mat($\cdot$): $
u = \text{vec}(U)$, $U = \text{mat}(u)$ where $u \in \mathbb{R}^{n_x  n_\xi}, U \in \mathbb{R}^{n_x \times n_\xi}$. 
A solution $u$ can be represented by a sum of vectors of tensor structure, or equivalently $U = \text{mat}(u)$  can be represented by a sum of rank-one matrices,
\begin{align}
u = \sum_{k=1}^{\kappa_u} z_k \otimes y_k  \: \Leftrightarrow \:U = \sum_{k=1}^{\kappa_u} y_k z_k^T = Y_{\kappa_u} Z_{\kappa_u}^T\label{sol_mat}
\end{align}
where $y_i \in \mathbb{R}^{n_x}$, $z_i \in \mathbb{R}^{n_\xi}$, and  $Y_{\kappa_u} =  
[y_1,\, \dots,\, y_{\kappa_u}] \in \mathbb{R}^{n_x \times \kappa_u}$ and $Z_{\kappa_u} = [z_1,\, \dots,\,z_{\kappa_u}]$
$\in \mathbb{R}^{n_\xi \times \kappa_u}$. 
If $\kappa_u$ is the rank of $U$, then we use $\kappa_u$ to refer to the rank of $u$ given in tensor structure; thus, $u$ is the sum of terms of rank-one tensor structure. With this notation, the stochastic Galerkin solution  $u_{hp}(\pmb{x},\,\xi)$ can be represented as 
\begin{align}
\label{sg_matrix_form}
\begin{split}
u_{hp}(\pmb{x},\,\xi)  = \Phi(\pmb{x})^T Y_{\kappa_u}Z_{\kappa_u}^T \Psi (\xi) =  \left( Y_{\kappa_u}^T \Phi(\pmb{x}) \right)^T  \left( Z_{\kappa_u}^T \Psi (\xi) \right)
\end{split}
\end{align}
where $\Phi(\pmb{x}) = [\phi_1(\pmb{x}),\, \dots,\, \phi_{n_x}(\pmb{x})]^T$ and $\Psi(\xi) = [\psi_1(\xi),\, \dots,\, \psi_{n_\xi}(\xi)]^T$. Note that as shown in \cite{tamellini2014model}, \eqref{sg_matrix_form} corresponds to a separated representation \cite{beylkin2005algorithms},
\begin{equation} 
u_{hp}(\pmb{x},\,\xi) = \sum_{i=1}^{\kappa_u} \hat{y}_i(\pmb{x}) \hat{z}_i (\xi),
\end{equation}  
where $\hat{y}_i(\pmb{x}) =  (\Phi(\pmb{x}))^T y_i$ and $\hat{z}_i(\xi) =  (\Psi(\xi))^T z_i$. We will use this representation to construct a new rank-reduction operator. In the discrete model \eqref{sg_matrix_form}, the rank of the solution is typically $\kappa_u = \min(n_x,\,n_\xi)$. 

In \cite{benner2013low, grasedyck2004existence}, it was shown that the solution to \eqref{linear_sys_tensor} can be approximated well by a quantity $\tilde{u}$ of rank $\kappa_{\tilde{u}} \ll \min(n_x,\,n_\xi)$  if the system matrix and the right-hand side has Kronecker-product structure. Thus, we seek a low-rank approximation to the solution $\tilde{u}$ to \eqref{linear_sys_tensor} for which
\begin{align}
A\tilde{u} = \left( \sum_{l=0}^M G_l \otimes K_l \right) \left( \sum_{k=1}^{\kappa_{\tilde{u}}} \tilde{z}_i \otimes \tilde{y}_i \right) \approx g_0 \otimes f_0  \label{sg_tensor}.
\end{align}

\subsection{Basic operations in tensor notation} 
We point out here a feature of the basic operations required by Krylov subspace methods in the setting we are considering, where the operators and data of interest have tensor format. The $m$th step of such methods results in the {\em Krylov subspace}, $\mathcal{K}_m(A,\,v_1) = \text{span}\{v_1,\,Av_1,\,\dots,\,A^{m-1}v_1\}$, which is  generated using matrix-vector products and addition/subtraction of vectors.

The matrix-vector product in \eqref{sg_tensor} can be represented as a sum of rank-one tensors by exploiting the properties of the Kronecker product,
\begin{equation}
A{u} = \sum_{l=0}^M \sum_{k=1}^{\kappa_{{u}}}  G_l z_k \otimes K_l y_k =\sum_{i=1}^{(M+1)\kappa_{{u}}} \hat{z}_i \otimes \hat{y}_i. 
\label{matvec_tensor}
\end{equation}
The latter expression in \eqref{matvec_tensor} suggests that in tensor notation, the matrix-vector product typically results in a vector with a higher rank.  Similarly, the addition of two vectors $u$ and $v$ of rank $\kappa_u$ and $\kappa_{{v}}$ in tensor notation gives
\begin{equation}
u + v = \sum_{i=1}^{\kappa_{{u}}} z_i \otimes y_i + \sum_{j= 1}^{\kappa_{{v}}} \hat{z}_j \otimes \hat{y}_j = \sum_{i=1}^{\kappa_{u} + \kappa_{v}} z_i \otimes y_i,
\end{equation}
where $y_{i+\kappa_u} = \hat{y}_i$ and $z_{i+\kappa_u} = \hat{z}_i$, $i=1,\,\dots,\,\kappa_v$, so that the resulting sum may have rank as large as $\kappa_u + \kappa_v$. Thus, although the goal is to find an approximate solution to \eqref{linear_sys_tensor} of low rank, two of the fundamental operations used in Krylov subspace methods tend to increase the rank of the quantities produced. Following \cite{ballani2013projection}, we will address this point in the next section.

\section{A preconditioned projection method in tensor format} \label{sec_lr_proj}
As is well known, the generalized minimum residual method (GMRES)  \cite{saad1986gmres} constructs an approximate solution $u_m \in u_0 + \mathcal{K}_m(A,\,v_1)$ where $u_0$ is an initial vector with residual $r_0 = f - Au_0$, $v_1 = r_0 / \| r_0 \|_2$, and $\mathcal{K}_m$ is the Krylov space. This is done by generating $V_m  = [v_1,\,\dots,\,v_m]$, where $\{v_j \}_{j=1}^m$ is an orthogonal basis for $\mathcal{K}_m$, and then computing $u_m$ whose residual $r_m$ is orthogonal to $W_m = AV_m$. The method is shown in Algorithm \ref{gmres_alg}. In this section, we discuss a variant of this method based on low-rank projection, where advantage is taken of the tensor format of the matrix $A$ and low-rank structure of the solution $u$.

\begin{algorithm}
\caption{GMRES method without restarting \cite{saad2003iterative}}
\label{gmres_alg}
\begin{algorithmic}[1]
\State set the initial solution $u_0$
\State $r_0 := f - Au_0$
\State $\tilde{v}_1 := r_0$
\State $v_1 := {\tilde{v}_1}/{\| \tilde{v}_1 \|}$
\For{$j = 1,\,\dots,\,m$}
\State $w_j := Av_j$
\State solve $(V_j^T V_j)\alpha = V_j^Tw_j$
\State $\tilde{v}_{j+1} :=  w_j - \sum_{i=1}^{j} \alpha_i v_i  $
\State $v_{j+1} := \tilde{v}_{j+1} / \| \tilde{v}_{j+1}\|$
\EndFor
\State solve $(W_m^T AV_m) y = W_m^T r_0$
\State $u_{m} := u_0 + V_m y$
\end{algorithmic}
\end{algorithm}

\subsection{Low-rank projection method with restarting} 
\label{low-rank-sec}
As we observed in Section \ref{sec_sg_model}, matrix-vector products and vector sums in tensor structure tend to increase the rank of the resulting objects. Thus, although we seek a solution of low rank, straightforward use of the GMRES method may lead to approximate solutions of higher rank than the desired solutions. This complication can be addressed using \textit{truncation operators} \cite{ballani2013projection, kressner2010krylov, kressner2011low, matthies2012solving, schwab2011sparse}, whereby vectors of high rank are replaced by ones of low rank. The truncation is inserted into the GMRES algorithm and is interleaved with the basic operations such as matrix-vector product and addition so that the ranks of the vectors used in the algorithm are kept low. 

\begin{algorithm}
\caption{Restarted low-rank projection method in tensor format}
\label{lrp_alg}
\begin{algorithmic}[1]
\State set the initial solution $\tilde{u}_0$
\For{ $k$ = $0,\,1,\,\dots$} 
\State $r_k := f - A\tilde{u}_k$
\If{$\| r_k \| / \| f \| < \epsilon$}
\State return $\tilde{u}_k$
\EndIf
\State $\tilde{v}_1 := \mathcal{T}_{\kappa}(r_k)$
\State $v_1 := {\tilde{v}_1}/{\| \tilde{v}_1 \|}$
\For{$j = 1,\,\dots,\,m$}
\State $w_j := Av_j$ \label{lrp-11}
\State solve $(V_j^T V_j)\alpha = V_j^Tw_j$ \label{lrp-12}
\State $\tilde{v}_{j+1} := \mathcal{T}_{\kappa}\left( w_j - \sum_{i=1}^{j} \alpha_i v_i \right) $ \label{lrp-13}
\State $v_{j+1} := \tilde{v}_{j+1} / \| \tilde{v}_{j+1}\|$ \label{lrp-14}
\EndFor
\State solve $(W_m^T AV_m) \beta = W_m^T r_k$ \label{lrp-16}
\State $\tilde{u}_{k+1} := \mathcal{T}_{\kappa} (\tilde{u}_k + V_m \beta)$
\EndFor
\end{algorithmic}
\end{algorithm}

Algorithm \ref{lrp_alg} summarizes the restarted low-rank projection method in tensor format \cite{ballani2013projection}.  As in the standard Arnoldi iteration used by GMRES, a new vector $\hat{v}_{j+1}$ is constructed by applying the linear operator $A$ to the previous basis vector $v_j$ and orthogonalizing the new basis vector $w_j$ with respect to the previous basis vectors $\{v_i\}_{i=1}^j$. The resulting vector is truncated to a vector  $\tilde{v}_{j+1}$ of low rank and normalized to $v_{j+1}$, which is then added to the set of basis vectors. The truncation operator $\mathcal{T}_{\kappa}$ truncates a tensor of higher rank to one of rank $\kappa$. Thus, all the basis vectors $\{v_i\}_{i=1}^{m}$ are of the same rank, $\kappa$. The basis vectors determine the subspace $\mathcal{K}_m = \text{span} \{v_1,\,\dots,\,v_{m} \}$, but because of truncation the basis vectors are not orthogonal and $\mathcal{K}_m$ is not a Krylov subspace. However, it is still possible to project the residual onto the subspace $\mathcal{W}_m = \text{span} \{w_1,\,\dots,\,w_m \}$ to find out whether the residual can be decreased by forming a new iterate $\tilde{u}_k + V_m \beta$. Note that all the vectors used in the entire iteration process are stored as the product of two matrices in the form like that shown in the right side of \eqref{sol_mat}. The ranks of these vectors will be discussed below.


\subsection{Preconditioned low-rank projection method}
To speed the convergence of the projection method, we consider a right-preconditioned system: 
\begin{equation}
AM^{-1}\hat{u} = f,\qquad \hat{u} = Mu.
\end{equation}
 For the stochastic diffusion problem, we consider $M = G_0 \otimes K_0$ as the precondtioner, which is known as the mean-based preconditioner \cite{powell2009block}. For the practical application of the preconditioner, we employ algebraic multigrid methods \cite{ruge1987algebraic}, where the action of $K_0^{-1}$ is replaced by $\tilde{K}_0^{-1}$, an application of a single V-cycle of an algebraic multigrid method. The preconditioned matrix-vector product is then
\begin{equation*}
 AM^{-1} \hat{u} = \sum_{l=0}^M \sum_{k=1}^{\kappa_{\hat{u}}}  G_l  \hat{z}_k \otimes  K_l \tilde{K}_0^{-1} \hat{y}_k, \quad \hat{u} = M{u} = \sum_{i=1}^{\kappa_{\hat{u}}} \hat{z}_i \otimes \hat{y}_i.
\end{equation*}
Note that $G_0^{-1}$ is the identity matrix because of the orthonormality of the stochastic basis functions. With right preconditioning and this preconditioner, the strategy for handling tensor rank is largely unaffected by preconditioning.


\section{Truncation methods} \label{sec_trunc}
As discussed in Section \ref{low-rank-sec}, in the low-rank projection method, truncation of tensors is  essential for the efficient computation of approximate solutions. In this section, we discuss the conventional approach for truncation and we introduce a new multilevel truncation method based on a coarse-grid solution.
\subsection{Truncation based on singular values}
Given a matricized vector $U = Y_{\kappa'} Z_{\kappa'}^T$  of rank $\kappa'$, a standard approach for truncation \cite{ballani2013projection,matthies2012solving} is to compute the singular value decomposition (SVD) of $U$ and compress $U$ into an approximation of desired rank $\kappa \ll \kappa'$. This can be done efficiently by computing QR factorizations of $Y_{\kappa'}$ and $Z_{\kappa'}$: 
\begin{equation*}
Y_{\kappa'} = Q_Y R_Y \in \mathbb{R}^{n_x \times \kappa'}, \quad Z_{\kappa'} = Q_Z R_Z \in \mathbb{R}^{n_\xi \times \kappa'}.
\end{equation*}
Then, one can compute the SVD of $R_Y R_Z^T$: 
\begin{equation*} 
R_Y R_Z^T = \hat{U}_{\kappa'}\hat{\Sigma}_{\kappa'}\hat{V}_{\kappa'}^T = \sum_{k=1}^{\kappa'} \hat{\sigma}_k \hat{u}_k \hat{v}_k^T
\end{equation*}
and truncate the sum with $\kappa$  terms to produce 
\begin{equation*}
\tilde{Y}_{\kappa} = Q_Y \hat{U}_{\kappa}  \hat{\Sigma}_{\kappa} \in \mathbb{R}^{n_x \times \kappa},\quad \tilde{Z}_{\kappa} = Q_Z \hat{V}_{\kappa} \in \mathbb{R}^{n_\xi \times \kappa}.
\end{equation*}
The truncated approximation of $U$ is then $\tilde{U} = \tilde{Y}_{\kappa} \tilde{Z}_{\kappa}^T$. The computational complexity of the truncation is $O( (n_x + n_\xi + \kappa) (\kappa')^2)$  \cite{golub2012matrix}, which grows quadradically with respect to $\kappa'$. In the next section, we introduce a new truncation method that avoids this computation.

\subsection{Truncation based on mutlilevel rank-reduction}
We now propose a multilevel rank-reduction strategy. We obtain insight into the rank structure of the solution using a coarse spatial grid computation. Then, we define a truncation operator based on the information obtained from this coarse-grid computation.

Let $u^c(\pmb{x},\,\xi) $ represent a solution obtained on a coarse spatial grid (i.e., $n_x$ is small). As in \eqref{sg_matrix_form}, $u^c(\pmb{x},\,\xi)$ can be represented as 
\begin{equation*}
\label{sg_coarse}
u^{c}(\pmb{x},\,\xi)  = \left(\Phi^c(\pmb{x})\right)^T U^c \Psi (\xi) = \left( (Y^c)^T \Phi^c(\pmb{x})\right)^T  \left( (Z^c)^T \Psi (\xi) \right).
\end{equation*}
Here, we propose to use $Z^c$ to define a truncation operator for use in the projection method to compute a solution for the problem on a finer grid. That is, the truncation operator is defined such that, given a matricized vector $U = Y_{\kappa'}Z_{\kappa'}^T$ of rank $\kappa'$, 
\begin{equation}
\label{trunc_multi}
\mathcal{T}_\kappa(U) \equiv \left( Y_{\kappa'}Z_{\kappa'}^T Z^c_\kappa \right) \left( Z^c_\kappa \right)^T = \tilde{U}
\end{equation}
where the resulting quantity $\tilde{U} = \tilde{Y}_\kappa\tilde{Z}_\kappa^T$ is of rank $\kappa$,
\begin{equation*}
\tilde{Y}_\kappa = Y_{\kappa'}Z_{\kappa'}^T Z^c_\kappa \in \mathbb{R}^{n_x \times \kappa}, \quad \tilde{Z}_\kappa = Z^c_\kappa \in \mathbb{R}^{n_\xi \times \kappa}.
\end{equation*}
The desired rank $\kappa$ is determined such that the relative residual $\|f^c - A^cu^{c,\, \kappa} \|_2 / \| f^c \|_2$ is smaller than a certain tolerance $\epsilon$ where $u^{c,\,\kappa}$ is a $\kappa$-term approximation of $u^c$. This truncation operation requires two matrix-matrix products, and the computational complexity of truncating a vector from $\kappa'$ to $\kappa$ is $O(\kappa' \kappa (n_x + n_\xi))$. 

For efficient coarse-grid computation, we use the Proper Generalized Decomposition (PGD) method developed in \cite{nouy2009generalized, tamellini2014model}, which computes a separated representation of a coarse-grid solution:  
\begin{equation}
\label{sep_rep}
u^{c,\,\kappa} (\pmb{x},\,\xi) = \sum_{i=1}^{\kappa} \tilde{y}_i (\pmb{x}) \tilde{z}_i (\xi).
\end{equation}
With the stochatic Galerkin discretization, each function can be represented as 
\begin{equation*}
\tilde{y}_i (\pmb{x}) = \sum_{k=1}^{n_x} \tilde{y}^{(i)}_k  \phi^c_k(\pmb{x}), \quad  \tilde{z}_i(\xi) = \sum_{l=1}^{n_\xi} \tilde{z}^{(i)}_l  \psi_l (\xi).
\end{equation*}
As a result, as in \eqref{sg_matrix_form}, 
\begin{equation*}
u^{c,\,\kappa} (\pmb{x},\,\xi) = \left( (\tilde{Y}_\kappa^c)^T \Phi^c(\pmb{x})\right)^T  \left( (\tilde{Z}_\kappa^c)^T \Psi (\xi) \right)
\end{equation*}
where $\tilde{Y}_{\kappa}^c = [\tilde{y}^{(1)},\, \cdots,\, \tilde{y}^{({\kappa})}] \in \mathbb{R}^{n_x \times \kappa}$ and $\tilde{Z}_{\kappa}^c = [\tilde{z}^{(1)},\, \cdots,\, \tilde{z}^{(\kappa)}] \in \mathbb{R}^{n_\xi \times \kappa}$ are coefficient matrices such that the $i$th elements of $\tilde{y}^{(j)}$ and $\tilde{z}^{(j)}$ are $\tilde{y}_i^{(j)}$ and $\tilde{z}_i^{(j)}$, respectively. Note that the new stochastic reduced basis $\hat{\psi}_j = (\tilde{z}^{(j)})^T  [ \psi_1(\xi),\,  \dots,\,  \psi_{n_\xi}(\xi) ]$ are also orthonormal with respect to a given probability measure $\rho$. Now, the discrete solution $U^c$ in \eqref{sg_coarse} is approximated by 
\begin{equation*}
U^{c,\,\kappa} = \tilde{Y}_\kappa^c (\tilde{Z}_\kappa^c)^T,
\end{equation*} 
and we can obtain $Z_\kappa^c$ by computing the SVD of $U^{c,\,\kappa}$, 
\begin{equation*}
U^{c,\,\kappa} = \hat{U}\hat{\Sigma}\hat{V}^T = Y^c_\kappa (Z^c_\kappa)^T.
\end{equation*} 
We  briefly explain how the PGD method computes a $\kappa$-term approximation in the next section.

\subsection{Proper Generalized Decomposition method}
The PGD method is a successive rank-1 approximation method. That is, the method incrementally identifies the function pairs $(\tilde{y}_i(\pmb{x}),\, \tilde{z}_i(\xi))$ of \eqref{sep_rep} one at a time. Once $i$ such pairs have been computed, the next pair $(\tilde{y}_{i+1},\, \tilde{z}_{i+1})$ is sought in $X_h \times S_M$ by imposing Galerkin orthogonality with respect to the tangent manifold of the set of rank-one elements at $\tilde{y}_{i+1} \tilde{z}_{i+1}$, which is $\{\tilde{y}_{i+1} \zeta + \upsilon\tilde{z}_{i+1};\,\upsilon\in X_h,\,\zeta \in S_M\}$: find $\tilde{y}_{i+1} \tilde{z}_{i+1}$ such that $\forall (\upsilon,\,\zeta) \in X_h \times S_M$
\begin{equation}
\label{PGD_tot}
\begin{split} 
\left \langle \int_D a(\pmb{x},\,\xi ) \nabla( u^{c,\,i} + \tilde{y}_{i+1}\tilde{z}_{i+1}) \cdot \nabla (\tilde{y}_{i+1}\zeta + \upsilon \tilde{z}_{i+1})  \right\rangle = \left \langle \int_D f  (\tilde{y}_{i+1}\zeta + \upsilon \tilde{z}_{i+1}) \right \rangle.
\end{split}
\end{equation}
It follows from \eqref{PGD_tot} that each component of a pair $(\tilde{y}_{i+1},\, \tilde{z}_{i+1})$ can be computed by solving two coupled problems: a deterministic problem \eqref{PGD_d} and a stochastic problem \eqref{PGD_s}. The deterministic problem is as follows: given $\tilde{z}_{i+1}$, find $\tilde{y}_{i+1}  \in X_h$ such that 
\begin{align}
\label{PGD_d}
\begin{split} 
\left \langle \int_D a(\pmb{x},\,\xi ) \nabla( u^{c,\,i} + \tilde{y}_{i+1}\tilde{z}_{i+1}) \cdot \nabla (\phi^c_j \tilde{z}_{i+1})  \right\rangle = \left \langle \int_D f \phi^c_j \tilde{z}_{i+1} \right \rangle,\quad j = 1,\,\dots,\,n_x^c.
\end{split}
\end{align}
The first basis function $\eta_1$ can be chosen arbitrarily at the beginning of the PGD method. The finite element discretization of $u_{i+1}$ yields a linear system of order $n_x^c$. Analogously, the stochastic problem starts with $\tilde{y}_{i+1}$ and finds $\tilde{z}_{i+1} \in S_M$ such that 
\begin{equation}
\label{PGD_s}
\left \langle \int_D a(\pmb{x},\,\xi ) \nabla( u^{c,\,i} + \tilde{y}_{i+1} \tilde{z}_{i+1}) \cdot \nabla (\tilde{y}_{i+1} \psi_{j})  \right\rangle = \left \langle \int_D f \tilde{y}_{i+1} \psi_{j} \right \rangle,\quad j = 1,\,\dots,\,n_\xi.
\end{equation}
Since $\tilde{z}_{i+1}$ is approximated by the PCE, $n_\xi$ unknowns have to be determined by solving a linear system of order $n_\xi$. 

Solutions of these sets of $\kappa$ systems of order $n_x^c$ and $\kappa$ systems of order $n_\xi$ produce   the $\kappa$-term approximation to the solution. The PGD method seeks solution pairs until the relative residual of the computed solution satisfies a given tolerance, 
\begin{equation*}
\|f^c - A^cu^{c,\, \kappa} \|_2 / \| f^c \|_2 < \epsilon.
\end{equation*} 
The accuracy of the $\kappa$-term approximation can also be improved by solving a set of $\kappa$ coupled equations: given $\{\tilde{y}_i \}_{i=1}^{\kappa}$, find $\{\tilde{z}_i\}_{i=1}^{\kappa}$ such that 
\begin{equation}
\label{PGD_u}
\left \langle \int_D a(\pmb{x},\,\xi ) \nabla( u^{(\kappa)} ) \cdot \nabla (\tilde{y}_i \psi_{j})  \right\rangle = \left \langle \int_D f \tilde{y}_i \psi_{j} \right \rangle,\quad i= 1,\,\dots,\,\kappa,\; j = 1,\,\dots,\,n_\xi. 
\end{equation}
This update requires the solution of a linear system of order $\kappa n_\xi$. For the stochastic diffusion problems, the update problem is solved once at the end of the PGD method. Note that the update problem could also be formulated for finding the deterministic components $\{u_i\}_{i=1}^{\kappa}$ if $n_x \ll n_\xi$, which requires a solution of a linear system of order $\kappa n_x$.



With the proposed truncation strategy, Algorithm \ref{overall_alg} summarizes the entire procedure to compute a solution on a finer grid. 

\begin{algorithm}
\caption{Preconditioned low-rank projection method with the multilevel rank-reduction}
\label{overall_alg}
\begin{algorithmic}[1]
\State Compute $u^{c,\,\kappa}$ that satisfies $\frac{\|f^c -  A^c u^{c,\,\kappa} \|_2}{\| f^c \|_2} < \epsilon$ using the PGD method
\State Compute $Z^c_\kappa$ such that $U^{c,\,\kappa} = Y^c_\kappa (Z^c_\kappa)^T$ and define $\mathcal{T}_\kappa(U) \equiv \left( U Z^c_\kappa \right) \left(Z^c_\kappa \right) ^T$
\State Run Algorithm \ref{lrp_alg} with $\mathcal{L} = AM^{-1}$, $f$, and $\mathcal{T}_\kappa$
\end{algorithmic}
\end{algorithm}

\section{Numerical experiments} \label{sec_num_exp}
In this section, we present the results of numerical experiments in which the proposed iterative solver is applied to some benchmark problems. The implementation of the spatial discretization is based on the Incompressible Flow and Iterative Solver Software (IFISS) package \cite{ifiss}. Example problems are posed on a square domain and $\ell$ is the spatial discretization parameter (i.e., $n_x = (2^\ell + 1)^2$).

For $a(\pmb{x},\xi)$ in \eqref{kl_exp}, we consider independent random variables $\{ \xi_i \}_{i=1}^M$ that are uniformly distributed over $[-\sqrt{3},\,\sqrt{3}]$, $a_0=1$ and unless otherwise specified, $\sigma=0.05$. As the covariance kernel, we use 
\begin{equation}
\label{exp_cov}
C(\pmb{x},\, \pmb{y}) =\sigma^2 \exp \left (  - \frac{|x_1 - y_1|}{c} - \frac{|x_2 - y_2|}{c} \right)
\end{equation}
where $c$ is the correlation length. The number of terms $M$ in the truncated expansion \eqref{kl_exp} is determined such that 95\% of the total variance is captured by $M$ terms (i.e., $( \sum_{i=1}^M \lambda_i ) / \left( \sum_{i=1}^{n_x} \lambda_i \right) > 0.95$). We use bilinear $Q_1$ elements to generate the finite element basis and Legendre polynomials as the stochastic basis functions because the underlying random variables have a uniform distribution. The default setting of the maximal polynomial degree $p$ is 3.


\subsection{Stochastic diffusion problem}
We consider the steady-state stochastic diffusion equation in \eqref{galerkin_form} on a domain $D = [0,\,1] \times [0,\,1]$ with forcing term $f(\pmb{x})=1$  and homogeneous Dirichlet boundary conditions, $u(\pmb{x},\,\omega)=0$ on $\partial D \times \Gamma$.

\vspace{0.5mm}
\textbf{Coarse spatial grid computation.} We compute $\kappa$-term approximations using the PGD method on a coarser spatial grid. Here $\ell^c$ is the refinement level for the coarse grid and $n^c_x$ is the number of degrees of freedom in the corresponding spatial domain excluding boundary nodes. 
We discuss choices of coarse spatial grid in Section \ref{sec_coarse_sel}. Table \ref{tab_coarse1} shows the rank $\kappa$ of solutions that satisfy the tolerance $\epsilon$ for varying correlation lengths $c$ and $M$ and the corresponding computation time $t_c$.

\begin{table}[htbp]
\caption{Rank ($\kappa$) of coarse-grid solutions satisfying a specified tolerance $\epsilon$ for the PGD computation, and CPU time ($t_c$) for coarse-grid computation using the PGD method, for varying $c$ and $M$.}
\begin{center}\footnotesize
\renewcommand{\arraystretch}{1.3}
\label{tab_coarse1}
\begin{tabular}{  c | c  c  c  c || c c c c }
\hline
 & \multicolumn{4}{c||}{$\epsilon = 10^{-5}$} &  \multicolumn{4}{c}{$\epsilon = 10^{-6}$}\\
\hline
 $c$ & 4 & 3 & 2.5 & 2 & 4 & 3 & 2.5 & 2 \\
 $M$,  $n_\xi$ &  5, 56 & 7, 120 & 10, 286 & 15, 816 & 5, 56 & 7, 120 & 10, 286 & 15, 816\\
\hline
 $n_x^c$$(\ell^c)$ & 225(4) & 225(4) & 961(5) & 961(5)& 225(4) & 225(4) & 961(5) & 961(5)\\
 Rank($\kappa$) & 25 & 40 & 65 & 115 & 35 & 65 & 100 & 210\\
CPU time($t_{c}$) & 2.49 & 3.47 & 8.35 & 45.08 & 2.93 & 5.04 & 14.83 & 162.71 \\
\hline
\end{tabular}
\end{center}
\end{table}

\textbf{Fine spatial grid computation.} With the truncation operator $\mathcal{T}_\kappa$ \eqref{trunc_multi} obtained from the coarse-grid solution (i.e., $Z^c_\kappa $), we solve the same stochastic diffusion problems on finer spatial grids $\ell = \{7,\,8,\,9\}$. For the fine-grid low-rank solutions, we use the rank $\kappa$ obtained from the coarse-grid solutions. For example, the third column of Table \ref{tab_fine1} shows the time required to find solutions of rank 25 satisfying the relative residual tolerance $10^{-5}$ when the number of terms in \eqref{kl_exp} is $M=5$. In Algorithm \ref{lrp_alg}, we set $m=8$ (like restarted GMRES(8)). In examining performance, we identify the number of cycles, $k$, performed for the outer for-loop in Algorithm \ref{lrp_alg}; this means that the number of matrix-vector products (i.e., the number of times line 10 is executed) is $mk$. Tables \ref{tab_fine1} and \ref{tab_fine2} show the number of cycles, $k$, and the computation time in seconds needed to compute approximate solutions with the tolerance $\epsilon = 10^{-5}$ and $10^{-6}$, respectively, (see line 4 of Algorithm \ref{lrp_alg}). Here, $t$ is the total time and $t_f$ excludes the time to compute the coarse-grid solution, $t_c$. The fine-grid computation time, $t_f$, consists of algorithm execution time and preconditioner set-up time, $t_{setup}$.
The execution times show ``textbook'' behavior, i.e., they grow linearly with the size of the spatial grid.\footnote{ An exception to this statement is when both $M$ and $n_x$ are large. For these cases, the problem does not fit into physical memory and memory swap-in/out time dominates the execution time.} If required memory for running Algorithm \ref{lrp_alg} exceeds the resources of our computing environment, solutions could not be computed and we denote these cases by OoM for ``Out-of-Memory''. Table \ref{tab_dof} shows the number of degrees of freedom of the fine spatial-grid problems for varying stochastic dimensions, $M$.

\begin{table}[!htbp]
\caption{CPU times to compute approximate solutions of the diffusion equation satisfying $\epsilon = 10^{-5}$  using the preconditioned low-rank projection method with the multilevel rank-reduction. Numbers of GMRES cycles are shown in parentheses.}
\begin{center}\footnotesize
\renewcommand{\arraystretch}{1.3}
\label{tab_fine1}
\begin{tabular}{ c | c | l l l l || c }
\hline
$n_x$$(\ell)$ &  & \multicolumn{1}{c}{$M$=5} & \multicolumn{1}{c}{$M$=7} & \multicolumn{1}{c}{$M$=10} & \multicolumn{1}{c||}{$M$=15} & $t_{setup}$\\
\hline
\multirow{2}{*}{\parbox{6mm}{\centering $129^2$ $(7)$ } } & $t_f$ &  \phantom{11}4.12 (1) &  \phantom{11}7.22 (1)&  \phantom{1}18.79 (1)  & \phantom{11}86.29 (1) & \phantom{1}1.76\\
& $t$ & \phantom{11}8.35  & \phantom{1}12.43 & \phantom{1}28.88 &  \phantom{1}132.15 \\
\hline
\multirow{2}{*}{\parbox{6mm}{\centering $257^2$ $(8)$ } }& $t_f$ & \phantom{1}12.55 (1) & \phantom{1}24.70 (1) & \phantom{1}74.71 (1) &  \phantom{1}330.45 (1)& 10.16\\
& $t$ & \phantom{1}25.17  & \phantom{1}38.37 & \phantom{1}93.20 &  \phantom{1}385.59\\
\hline
\multirow{2}{*}{\parbox{6mm}{\centering $513^2$ $(9)$ } }& $t_f$ & \phantom{1}92.83 (1)&  102.42 (1)& 353.07 (1) & 2717.03 (1) & 92.41\\
& $t$ & 147.17  & 197.87 & 453.71 &  2854.62 \\
\hline
\end{tabular}
\end{center}
\end{table}

\begin{table}[htbp]
\caption{CPU times to compute approximate solutions of the diffusion equation satisfying $\epsilon = 10^{-6}$  using the preconditioned low-rank projection method with the multilevel rank-reduction. Numbers of GMRES cycles are shown in parentheses.}
\begin{center}\footnotesize
\renewcommand{\arraystretch}{1.3}
\label{tab_fine2}
\begin{tabular}{ c | c | l l l l  || c}
\hline
$n_x$$(\ell)$ &  & \multicolumn{1}{c}{$M$=5} & \multicolumn{1}{c}{$M$=7} & \multicolumn{1}{c}{$M$=10} & \multicolumn{1}{c}{$M$=15} & $t_{setup}$\\
\hline
\multirow{2}{*}{\parbox{6mm}{\centering $129^2$ $(7)$ } }& $t_f$ &  \phantom{11}5.40 (1) & \phantom{1}12.50 (1) & \phantom{11}35.09 (1) & \phantom{1}233.54 (1) & \phantom{1}1.79\\
& $t$ &  \phantom{1}10.14  & \phantom{1}19.32 & \phantom{11}51.69 & \phantom{1}398.06\\
\hline
\multirow{2}{*}{\parbox{6mm}{\centering $257^2$ $(8)$ } }& $t_f$ &  \phantom{1}17.23 (1) & \phantom{1}46.07 (1) & \phantom{1}137.19 (1) & 1004.40 (1) & 10.53\\
& $t$ &  \phantom{1}30.55  & \phantom{1}61.41 & \phantom{1}162.90& 1177.68\\
\hline
\multirow{2}{*}{\parbox{6mm}{\centering $513^2$ $(9)$ } }& $t_f$ & \phantom{1}70.37 (1)& 217.12 (1) & 1225.77 (1) & \multicolumn{1}{c||}{OoM} & 92.81\\
& $t$ & 166.24  & 315.18 & 1333.63 & \multicolumn{1}{c||}{OoM}\\
\hline
\end{tabular}
\end{center}
\end{table}

\begin{table}[!htbp]
\caption{Number of degrees of freedom of the fine-grid discretizations with $p=3$, for varying spatial-grid refinement level, $\ell$, and number of random variables, $M$.}
\begin{center}\footnotesize
\renewcommand{\arraystretch}{1.3}
\label{tab_dof}
\begin{tabular}{ c |  rrrr  }
\hline
$\ell$   & \multicolumn{1}{c}{$M$=5} & \multicolumn{1}{c}{$M$=7} & \multicolumn{1}{c}{$M$=10} & \multicolumn{1}{c}{$M$=15} \\
\hline
7  & 931,896 & 1,996,920 & 4,759,326  & 13,579,056 \\
\hline
8 & 3,698,744 &7,925,880 & 18,890,014 & 53,895,984\\
\hline
9 &14,737,464 & 31,580,280& 75,266,334& 214,745,904\\
\hline
\end{tabular}
\end{center}
\end{table}

\textbf{Example problems with varying $\sigma$ and $p$.} We examine the rank structure of the numerical solutions of the stochastic diffusion problems and assess the performance of the proposed solution algorithm for different values of maximal degree of stochastic polynomial, $p$ in \eqref{td_set}, and variance $\sigma^2$ of the random field $a(\pmb{x},\,\xi)$. As in the previous numerical experiments, we first identify the rank structure and define the truncation operator from coarse-grid computation. Then, we solve the same problems on a finer grid by using the proposed low-rank projection method with the multilevel rank-reduction scheme.

Table \ref{tab_degree} shows the computation time needed to compute approximate solutions of the stochastic diffusion problems with $M=7$ for varying maximal polynomial degree $p$. The required ranks of the approximate solutions are not affected by the number of terms in the polynomial expansion. However, the computation time is increased for the polynomial expansion with higher maximal polynomial degree because the size of $\{ G_i \}_{i=0}^{M}$ and the size of the stochastic part of the solution gets larger as the number of terms in the PCE is increased. 
\begin{table}[htbp]
\caption{CPU times $t$ to compute approximate solutions of the diffusion equation satisfying $\epsilon = 10^{-5}$ and $10^{-6}$ using the preconditioned low-rank projection method with the multilevel rank-reduction for varying maximal polynomial degree $p$ (stochastic degrees-of-freedom, $n_\xi$, in the parenthesis).}
\begin{center}\footnotesize
\renewcommand{\arraystretch}{1.3}
\label{tab_degree}
\begin{tabular}{ c |  rrr  || rrr}
\hline
&  \multicolumn{3}{c||}{$\epsilon = 10^{-5}\; (\kappa=40)$} &\multicolumn{3}{c}{$\epsilon = 10^{-6}\; (\kappa=65)$}\\
\hline
\multirow{1}{*}{\parbox{10mm}{\centering $n_x$$(\ell)$ }}&   \multicolumn{1}{c}{$p = 3$ (120)} & $p = 4$ (330) & $p =5$ (792) & $p=3$ (120)& $p=4$ (330) & $p=5$ (792)\\
\hline
$129^2$$(7)$ &  12.43  & 15.55 & 21.56 & 19.32 &23.42 & 38.49\\
$257^2$$(8)$ &  38.37  & 44.27 & 56.79 & 61.41 & 69.17 & 91.10 \\
$513^2$$(9)$ &  197.87 & 217.38 & 252.39 & 315.18 &322.86 & 383.89\\
\hline
\end{tabular}
\end{center}
\end{table}

Table \ref{tab_sigma1} shows the computation time $t$ needed to compute approximate solutions of the stochastic diffusion problems that satisfy the tolerance $10^{-5}$ and $10^{-6}$ for varying variance, $\sigma^2$. In general, the example problem with a larger variance requires a higher rank to satisfy the stopping tolerance, which, therefore, requires more computational effort.


\begin{table}[htbp]
\caption{CPU times $t$ and rank $\kappa$ to compute approximate solutions of the diffusion equation satisfying $\epsilon = 10^{-5}$ and $10^{-6}$ using the preconditioned low-rank projection method with the multilevel rank-reduction, for varying $\sigma$.} 
\begin{center}\footnotesize
\renewcommand{\arraystretch}{1.3}
\label{tab_sigma1}
\begin{tabular}{ c |  c | rrrr  || rrrr }
\hline
& &  \multicolumn{4}{c||}{$\epsilon = 10^{-5}$} &\multicolumn{4}{c}{$\epsilon = 10^{-6}$}\\
\hline
$\sigma$ & $n_x$ & \multicolumn{1}{c}{$M$=5} & \multicolumn{1}{c}{$M$=7} & \multicolumn{1}{c}{$M$=10} & \multicolumn{1}{c||}{$M$=15} & \multicolumn{1}{c}{$M$=5} & \multicolumn{1}{c}{$M$=7} & \multicolumn{1}{c}{$M$=10} & \multicolumn{1}{c}{$M$=15}\\
\hline
\multirow{4}{*}{0.01}  & & $\kappa=15$ & $\kappa=20$ & $\kappa=35$ & $\kappa=55$  &$\kappa= 20$ & $\kappa=30$ & $\kappa=50$ & $\kappa=85$   \\
\cline{3-10}
& {$129^2$}  & 7.28 & 8.65  & 15.01  & 45.69 & 7.87 & 10.81  & 20.76 & 83.07   \\
&{$257^2$}& 21.47  & 26.08  & 47.21  & 135.75 & 23.30  & 31.94  & 66.92 & 240.98  \\
&{$513^2$}&  130.93  & 150.85  & 236.34  & 922.87 & 137.98  & 173.03 & 333.70  & 1893.89  \\
\hline
\multirow{4}{*}{0.05}  &  & $\kappa=25$ & $\kappa=40$ &$\kappa= 65$ & $\kappa=115$&$\kappa=35$ & $\kappa=65$ & $\kappa=100$ & $\kappa=210$ \\
\cline{3-10}
&{$129^2$} &  8.35   & 12.43  & 28.88  & 132.15 &  10.14  & 19.32  & 51.69  & 398.06 \\
&{$257^2$} & 25.17  & 38.37  & 93.20  & 385.59 &  30.55  & 61.41  & 162.90  & 1177.68 \\
&{$513^2$} &  147.17  & 197.87  & 453.71  & 2854.62 &  166.24   & 315.18  & 1333.63  & OoM\\
\hline
\multirow{4}{*}{0.1}  & & $\kappa=35$ &$\kappa= 60$ & $\kappa=100$ & $\kappa=180$ & $\kappa=50$ & $\kappa=85$ & $\kappa=145$ & -\phantom{1} \\ 
\cline{3-10}
& {$129^2$} & 9.78   & 17.24  & 50.70 & 297.35 & 8.79   & 28.37  & 113.53  & OoM\\
& {$257^2$} & 29.98   & 54.94  & 157.76 & 866.41 & 41.69   & 94.48  & 356.50  & OoM\\
& {$513^2$} & 164.48   & 273.33  & 1324.47 & OoM & 208.15   & 515.29 & 2902.95  & OoM\\
\hline
\end{tabular}
\end{center}
\end{table}

\textbf{Comparison to a truncation operator based on singular values.} We compare the performance of the proposed solver to the preconditioned low-rank projection method combined with the conventional truncation operator from \cite{kressner2011low}. Table \ref{tab_fine_svd_multi} shows the computation time required to compute approximate solutions using the conventional and new truncation strategies. The total computation time, $t$, of the low-rank projection method with the multilevel rank reduction includes both coarse-grid, $t_c$, and fine-grid computations, $t_f$. The low-rank projection method with the SVD-based truncation operator, which is implemented based on \cite{ballani2013projection}, does not  require a coarse-grid computation and can start with any arbitrary initial guess for rank, $\kappa$. For these computations, we used the values of rank identified in the coarse-grid computations, which are illustrated in Table \ref{tab_coarse1}, for the initial rank.

\begin{table}[htbp]
\caption{CPU times to compute approximate solutions of the diffusion equation satisfying $\epsilon = 10^{-5}$ and $10^{-6}$ using the preconditioned low-rank projection (LRP) methods with the multilevel rank-reduction and the singular value based truncation on the level 8 spatial grid (i.e., $n_x = 257^2$).}
\begin{center}\footnotesize
\renewcommand{\arraystretch}{1.3}
\label{tab_fine_svd_multi}
\begin{tabular}{c | c | c | r r r r r}
\hline
&Solver &  & \multicolumn{1}{c}{$M$=5} & \multicolumn{1}{c}{$M$=7} & \multicolumn{1}{c}{$M$=10} & \multicolumn{1}{c}{$M$=15} & \multicolumn{1}{c}{$M=20$} \\
 \hline
\multirow{2}{*}{$\epsilon = 10^{-5}$}&LRP-SVD& $t_{\text{SVD}}$ & 55.04  & 108.11 & 284.27 & 1280.65 & 5691.19\\
\cline{2-8}
&LRP-Multilevel&$t$ & 25.17  & 38.37 & 93.20 & 385.59 & 1943.49\\
\hline
\multirow{2}{*}{$\epsilon = 10^{-6}$}&LRP-SVD& $t_{\text{SVD}}$ & 76.03  & 198.20 & 564.12 & 5131.32 & OoM\\
\cline{2-8}
&LRP-Multilevel&$t$ & 30.55  & 61.41 & 162.90 & 1177.68 & OoM\\
\hline
\end{tabular}
\end{center}
\end{table}

\textbf{PGD as a solver on a finer spatial grid.} The PGD method could be applied directly to the fine-grid problems. We assess the performance of the PGD method for computing fine-grid solutions in Table \ref{tab_pgd_finer1}, which shows the rank and computation time for computing approximate solutions that satisfy the tolerance $10^{-5}$ and $10^{-6}$ using PGD on a finer spatial grid. For the low-rank projection method, we record total computation time, $t$, which includes coarse-grid computation, $t_c$, AMG preconditioner set-up, $t_{setup}$, and fine-grid computation time, $t_f$. We compare the rank and the computation time for computing solutions using the PGD method and the proposed multilevel projection method. 
The proposed low-rank projection method runs faster and requires somewhat smaller ranks than the PGD method.

\vspace{0.5mm}
\textbf{Remark.} We also tested the techniques compared in Tables \ref{tab_fine_svd_multi} and \ref{tab_pgd_finer1} for different values of $\sigma$, $\sigma=0.01$ and $0.1$, with similar results. Indeed, the performance of LRP-Multilevel is more favorable for the larger value $\sigma=0.1$.
\begin{table}[htbp]
\caption{CPU times to compute approximate solutions  of the diffusion equation satisfying $\epsilon=10^{-5}$ and $10^{-6}$ using the PGD method and the preconditioned low-rank projection methods with the multilevel rank-reduction on the level 8 spatial grid (i.e., $n_x = 257^2$).}
\begin{center}\footnotesize
\renewcommand{\arraystretch}{1.3}
\label{tab_pgd_finer1}
\begin{tabular}{c | c | c | c c c c c}
\hline
&Solver &  & $M=5$ & $M=7$ & $M=10$ & $M=15$ & $M=20$\\
\hline
\multirow{4}{*}{$\epsilon=10^{-5}$} &\multirow{2}{*}{PGD} &$\kappa$ & 25 & 45 & 65 & 125 & 195\\
&&$t$&  43.78 & 109.72 & 228.73 & 940.69 & 3066.87\\
\cline{2-8}
&\multirow{2}{*}{LRP-Multilevel} & $\kappa$ & 25 & 40 & 65 & 115 & 180\\
&& $t$ & 25.17  & 38.37 & 93.20 & 385.59 & 1943.49 \\
\hline
\multirow{4}{*}{$\epsilon=10^{-6}$} &\multirow{2}{*}{PGD} &$\kappa$ & 40 & 70 & 110 & 225 & OoM\\
&&$t$&  74.43 & 214.82 & 533.10 & 2713.70 & OoM\\
\cline{2-8}
&\multirow{2}{*}{LRP-Multilevel} & $\kappa$ & 35 & 65 & 100 & 210 & OoM\\
&& $t$ & 30.55  & 61.41 & 162.90 & 1177.68 & OoM \\
\hline
\end{tabular}
\end{center}
\end{table}

\subsection{Stochastic convection-diffusion problem}
\begin{figure}[htbp]
	\centering
	\includegraphics[angle=0, scale=0.45]{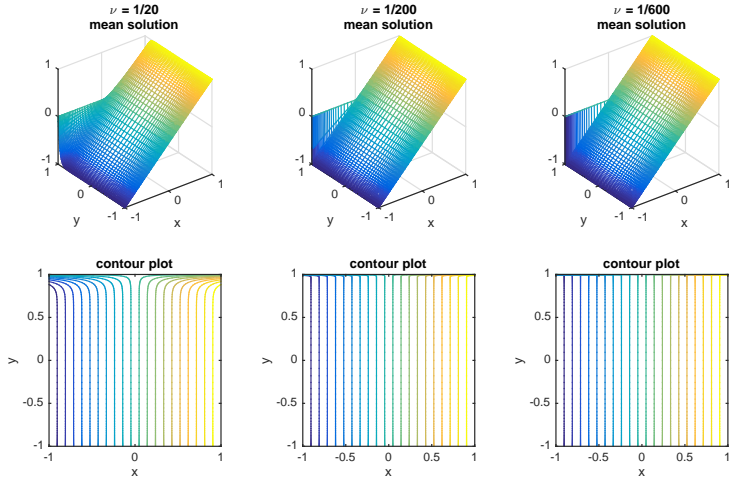}
	\caption{Mean solutions and their contour plots computed on the level 6 spatial grid for varying $\nu$.}
	\label{mean_sol}
\end{figure}

For a second benchmark problem, we consider the steady-state convection-diffusion equation defined on $D=[-1,\,1]\times[-1,\,1]$ with non-homogeneous Dirichlet boundary conditions, constant vertical wind $\vec{w} = (0,\,1)$, and $f=0$,
\begin{equation}
\label{cd_strong}
\left\{
\begin{array}{r l l}
\nu \nabla \cdot (a(\pmb{x},\,\xi) \nabla u(\pmb{x},\,\xi) ) + \vec{w} \cdot \nabla u(\pmb{x},\,\xi) &= f(\pmb{x},\,\xi) &\text{ in } D \times \Gamma,\vspace{2mm}\\
u(\pmb{x},\,\xi) &= g_D(\pmb{x}) &\text{ on } \partial D \times \Gamma,
\end{array}
\right.
\end{equation}
where $g_D(\pmb{x})$ is determined by
\begin{equation}
g_D(\pmb{x}) = 
\left\{
\begin{array}{c l c}
g_D(x,\,-1) = x, &g_D(x,\,1) = 0,\\
g_D(-1,\,y) = -1, & g_D(1,\,y) = 1,
\end{array}
\right.
\end{equation}
where the latter two approximations hold except near $y=1$, and $\nu$ is the viscosity parameter. We consider the convection-dominated case (i.e., $\nu < 1$) and employ the streamline-diffusion method for stabilization \cite{brooks1982streamline}. Here, we define the element \textit{Peclet} number 
\begin{equation}
\mathcal{P}_k = \frac{\|\vec{w}_k\|_2 h_k}{2\nu}
\end{equation}
where $\|\vec{w}_k\|_2$ is the $\ell_2$ norm of the wind at the element centroid and $h_k$ is a measure of the element length in the direction of the wind. Note that the solution has an \textit{exponential boundary layer} near $y=1$ where the value of the solution dramatically changes essentially from $-1$ to 0 on the left and $+1$ to 0 on the right \cite{elman2014finite}. Figure \ref{mean_sol} illustrates the mean of solutions $\langle u(\pmb{x},\xi) \rangle_\rho$ computed on the level 6 spatial grid and corresponding contour plots for varying viscosity parameter, $\nu$.

Given $a(\pmb{x},\,\xi)$ in \eqref{kl_exp}, we again discretize \eqref{cd_strong} using the finite element  method and the generalized PCE.  The result is a linear system in tensor product notation
\begin{align}
\label{linear_sys_tensor_cd}
\left(G_0 \otimes \nu K_0 + \sum_{l=1}^M G_l \otimes \nu  K_l + G_0 \otimes N + G_0 \otimes S \right) u = g_0 \otimes f_0
\end{align}
where the convection term $N$ and the streamline-diffusion term $S$ are given by
\begin{align*}
[N]_{ij} &= \int_D \vec{w} \cdot \nabla \phi_i (\pmb{x}) \phi_j (\pmb{x}) d\pmb{x}, \\
[S]_{ij} &= \sum_{k=1}^{n_e} \delta_k \int_D (\vec{w} \cdot \nabla \phi_i) (\vec{w} \cdot \nabla \phi_j) d,\pmb{x}
\end{align*}
$n_e$ is the number of elements in the finite element discretization, and 
\begin{equation}
\delta_k = 
\left\{
\begin{array}{c l c}
\frac{h_k}{2 \|\vec{w} \|_2} \left(1 - \frac{1}{\mathcal{P}_k} \right) & \text{if } \mathcal{P}_k > 1\\
0 &\text{if } \mathcal{P}_k \leq 1
\end{array}
\right..
\end{equation}
As the preconditioner, we choose $M \approx G_0 \otimes (K_0 + N + S)$ where the action of $(K_0 + N + S)^{-1}$ is replaced by application of a single V-cycle of an AMG method. In the PGD method, the non-homogeneous Dirichlet boundary condition is handled by introducing an extended affine space \cite{nouy2009generalized}: $u^{c} \approx u_{bc} + u^{c,\,\kappa}$ where $u_{bc}$ is the boundary nodal functions such as $u_{bc} = \sum_{k \in \partial D} u^{(bc)}_k \phi_k(\pmb{x})$. For the stochastic convection-diffusion problems, the update problems \eqref{PGD_u} need to be solved more often to compute an approximate solution of a desired accuracy with fewer terms.

\vspace{0.5mm}
\textbf{Numerical results.}  To cope with the existence of the exponential boundary layer in the solution, we use vertically stretched spatial grids. We examine the performance of the low-rank projection method for varying viscosity parameter $\nu$, and we set $m=10$ for Algorithm \ref{lrp_alg}. Table \ref{tab_cd1} and \ref{tab_cd2} show $\kappa$ computed by the PGD method, coarse-grid computation time $t_c$, and fine grid computation time $t_f$ to compute approximate solutions on fine spatial grids $\ell = \{7,\,8,\,9\}$ satisfying $10^{-5}$ and $10^{-6}$, respectively. Underlined numbers in the spatial grid level indicates cases where streamline diffusion is not needed.

\begin{table}[htbp]
\caption{CPU times to compute approximate solutions of the convection-diffusion equation satisfying $\epsilon = 10^{-5}$ using the preconditioned low-rank projection methods with the multilevel rank-reduction method for varying $\nu$. Numbers of GMRES cycles are shown in parentheses.} 
\begin{center}\footnotesize
\renewcommand{\arraystretch}{1.3}
\label{tab_cd1}
\begin{tabular}{ c | c | c | llll || r}
\hline
$\nu$ & $\ell$ &  & \multicolumn{1}{c}{$M=5$} &\multicolumn{1}{c}{$M= 7$} & \multicolumn{1}{c}{$M=10$} & \multicolumn{1}{c||}{$M=15$} & $t_{setup}$\\
\hline
\multirow{4}{*}{$\frac{1}{20}$}  &  \multirow{2}{*}{$4$} &$\kappa$ & \multicolumn{1}{c}{25} & \multicolumn{1}{c}{35} & \multicolumn{1}{c}{55} & \multicolumn{1}{c||}{65$^{\ast}$} \\
 & & $t_c$ & \phantom{11}2.56 &	\phantom{11}4.83 & 	\phantom{1}26.34 & \phantom{11}58.92$^{\ast}$ \\
 \cline{2-8}
&{$\underline{7}$} & $t_f$ & \phantom{11}5.73 (1)  & \phantom{11}9.47 (1)  & \phantom{1}24.86 (1)  & \phantom{11}72.29 (1) & 6.14\\
&{$\underline{8}$} & $t_f$ & \phantom{1}20.52 (1) &	\phantom{1}36.66 (1)  & \phantom{1}98.72 (1)& \phantom{1}248.31 (1) & 30.57\\
&{$\underline{9}$} & $t_f$ & \phantom{1}84.55 (1)&	 152.69 (1)& 592.63 (1) & 1953.52 (2) & 338.28\\
\hline
\multirow{5}{*}{$\frac{1}{100}$}  &  \multirow{2}{*}{$4$} &$\kappa$ & \multicolumn{1}{c}{20} & \multicolumn{1}{c}{25} & \multicolumn{1}{c}{45} & \multicolumn{1}{c||}{55$^{\ast}$} \\
 & & $t_c$ & \phantom{11}2.94 &	\phantom{11}3.12 &	\phantom{1}16.28 & \phantom{11}47.24$^{\ast}$\\
\cline{2-8}
&{$7$} & $t_f$ & \phantom{11}5.06 (1) &	 \phantom{11}7.28 (1) &	\phantom{1}18.90 (1)& \phantom{11}60.66 (1) & 6.34\\
&{$8$} & $t_f$ & \phantom{1}16.87 (1)  &	\phantom{1}26.36 (1)  &	\phantom{1}74.26 (1) & \phantom{1}202.29 (1) & 35.52\\
&{$\underline{9}$} & $t_f$ & 121.98 (2)  & 201.62 (2) & 745.92 (2)  & 3079.24 (2)& 341.41\\
\hline
\multirow{4}{*}{$\frac{1}{200}$}  &  \multirow{2}{*}{$5$} &$\kappa$ & \multicolumn{1}{c}{20} & \multicolumn{1}{c}{25} & \multicolumn{1}{c}{45} & \multicolumn{1}{c||}{50}\\
 & & $t_c$ & \phantom{11}2.91 &	\phantom{11}4.79 &	\phantom{1}16.54 & \phantom{11}46.85\\
\cline{2-8}
& {$7$} &$t_f$ & \phantom{11}5.16 (1)  & 	\phantom{11}7.21 (1) 	& \phantom{1}16.57 (1)  & \phantom{11}53.97 (1) & 6.35\\
& {$8$} &$t_f$ & \phantom{1}17.57 (1)  &	\phantom{1}25.05 (1)  & \phantom{1}63.56 (1)  & \phantom{1}175.30 (1) & 35.89\\
& {$9$} &$t_f$ & 123.73 (2)  &	200.10 (2) & 605.50 (2)  & 2568.41 (2) & 344.87\\
\hline
\multirow{5}{*}{$\frac{1}{400}$}  &  \multirow{2}{*}{$5$} &$\kappa$ & \multicolumn{1}{c}{20} & \multicolumn{1}{c}{20} & \multicolumn{1}{c}{35}  & \multicolumn{1}{c||}{45$^{\dagger}$} \\
 & & $t_c$ & \phantom{11}2.94 &	\phantom{11}3.79 &	\phantom{1}12.49 & \phantom{11}82.06$^{\dagger}$\\
\cline{2-8}
& {$7$} &$t_f$ & \phantom{11}8.61 (2)  &	 \phantom{11}9.84 (2) 	& \phantom{1}26.97 (2) & \phantom{11}85.01 (2) & 6.09 \\
& {$8$} &$t_f$ & \phantom{1}31.55 (2)  & \phantom{1}37.74 (2) 	& 111.31 (2) & \phantom{1}298.49 (2) & 34.93\\
& {$9$} &$t_f$ & 133.45 (2)&	158.01 (2)  &	512.88 (2) & 2080.60 (2) & 342.12\\
\hline
\multirow{4}{*}{$\frac{1}{600}$}  &  \multirow{2}{*}{$6$} &$\kappa$ & \multicolumn{1}{c}{20} & \multicolumn{1}{c}{20} & \multicolumn{1}{c}{35}  & \multicolumn{1}{c||}{45}\\
 & & $t_c$ & \phantom{11}9.79 &	\phantom{1}13.20 &	\phantom{1}34.47 & \phantom{11}94.79\\
\cline{2-8}
& {$7$} &$t_f$ &  \phantom{11}8.27 (2) & \phantom{1}10.07 (2) 	&  \phantom{1}26.91 (2) &  \phantom{11}82.30 (2) & 6.14 \\
& {$8$} &$t_f$ & \phantom{1}31.94 (2) &	\phantom{1}39.84 (2)	& 109.25 (2) & \phantom{1}295.25 (2) & 33.25 \\
& {$9$} &$t_f$ & 343.80 (2)& 163.90 (2) & 506.42 (2) & 1977.83 (2) & 342.98\\
\hline
\end{tabular}
\end{center}
\end{table}

\begin{table}[htbp]
\caption{CPU times to compute approximate solutions of the convection-diffusion equation satisfying $\epsilon = 10^{-6}$ using the preconditioned low-rank projection methods with the multilevel rank-reduction method for varying $\nu$. Numbers of GMRES cycles are shown in parentheses.} 
\begin{center}\footnotesize
\renewcommand{\arraystretch}{1.3}
\label{tab_cd2}
\begin{tabular}{ c | c | c | llll || r }
\hline
$\nu$ & $\ell$ &  & \multicolumn{1}{c}{$M=5$} &\multicolumn{1}{c}{$M= 7$} & \multicolumn{1}{c}{$M=10$} & \multicolumn{1}{c||}{$M=15$} & $t_{setup}$\\
\hline
\multirow{4}{*}{$\frac{1}{20}$}  &  \multirow{2}{*}{$4$} &$\kappa$ &  \multicolumn{1}{c}{35}	 & \multicolumn{1}{c}{50} &	 \multicolumn{1}{c}{75}   &  \multicolumn{1}{c||}{105$^{\ast}$} \\
 & & $t_c$ &  \phantom{11}3.31 &	\phantom{11}9.17 &	\phantom{11}60.51 & \phantom{11}194.33$^{\ast}$ \\
\cline{2-8}
&{$\underline{7}$} &$t_f$  & \phantom{1}13.92 (2) & \phantom{1}27.47 (2)  &	\phantom{11}80.78 (2)  & \phantom{11}275.96 (2) & 6.14\\
&{$\underline{8}$}&$t_f$   & \phantom{1}52.45 (2) & 106.11 (2) &	 \phantom{1}311.59 (2) & \phantom{1}1042.40 (2) & 30.57\\
&{$\underline{9}$}& $t_f$  & 220.67 (2)  & 534.61 (2)  &	2694.26 (2) & \phantom{1}8101.20 (2) & 338.28\\
\hline
\multirow{5}{*}{$\frac{1}{100}$}  &  \multirow{2}{*}{$4$} &$\kappa$   &  \multicolumn{1}{c}{30} &  \multicolumn{1}{c}{40} &  \multicolumn{1}{c}{65} & \multicolumn{1}{c||}{95$^{\ast}$}  \\
 & & $t_c$ & \phantom{11}2.83 &	\phantom{11}6.25 & \phantom{11}38.39 & \phantom{11}155.83$^{\ast}$ \\
\cline{2-8}
&{$7$} & $t_f$   & \phantom{1}12.34 (2)  & \phantom{1}21.28 (2)  &	\phantom{11}65.02 (2) & \phantom{11}239.91 (2) & 6.34\\
&{$8$} & $t_f$   & \phantom{1}46.67 (2)  & \phantom{1}85.66 (2)  & \phantom{1}255.79 (2) & \phantom{11}895.81 (2) & 35.52\\
&{$\underline{9}$} & $t_f$   & 273.45 (3) & 549.82 (3) & 3069.96 (3) & 10963.03 (3) & 341.41\\
\hline
\multirow{4}{*}{$\frac{1}{200}$}  &  \multirow{2}{*}{$5$} &$\kappa$ &  \multicolumn{1}{c}{25} & \multicolumn{1}{c}{40} & \multicolumn{1}{c}{60} &  \multicolumn{1}{c||}{85} \\
 & & $t_c$  & \phantom{11}3.46 & \phantom{11}8.57 & \phantom{11}38.35 & \phantom{11}122.49\\
\cline{2-8}
& {$7$} &$t_f$   & \phantom{1}10.52 (2)  & \phantom{1}21.43 (2)  &	\phantom{11}56.36 (2)  & \phantom{11}204.09 (2) & 6.35\\
& {$8$} &$t_f$   & \phantom{1}39.39 (2) &	 \phantom{1}84.14 (2) &	\phantom{1}219.36 (2) & \phantom{11}732.88 (2) & 35.89\\
& {$9$} &$t_f$   & 226.83 (3)  &	 547.62 (3)	 &  2627.98 (3) & \phantom{1}9284.60 (3)& 344.87\\
\hline
\multirow{5}{*}{$\frac{1}{400}$}  &  \multirow{2}{*}{$5$} &$\kappa$ &  \multicolumn{1}{c}{25} &  \multicolumn{1}{c}{35} &  \multicolumn{1}{c}{55} &  \multicolumn{1}{c||}{75$^{\dagger}$}  \\
 & & $t_c$ &  \phantom{11}3.49 &\phantom{11}6.63 &	\phantom{11}30.50 & \phantom{11}151.46$^{\dagger}$ \\
\cline{2-8}
& {$7$} &$t_f$   & \phantom{1}10.44 (2) &	 \phantom{1}17.96 (2)  & \phantom{11}50.96 (2) & \phantom{11}161.58 (2) & 6.09 \\
& {$8$} &$t_f$  & \phantom{1}40.02 (2)  &	 \phantom{1}70.82 (2) &	\phantom{1}204.71 (2) & \phantom{11}610.23 (2) & 34.93\\
& {$9$} &$t_f$  & 239.04 (3)  &	441.73 (3) &  2106.30 (3) & \phantom{1}7817.82 (3) & 342.12\\
\hline
\multirow{4}{*}{$\frac{1}{600}$}  &  \multirow{2}{*}{$6$} &$\kappa$ &  \multicolumn{1}{c}{30} &  \multicolumn{1}{c}{35} &  \multicolumn{1}{c}{45}&  \multicolumn{1}{c||}{65}\\
 & & $t_c$ & \phantom{1}17.99 &	 \phantom{1}22.03 &	\phantom{11}47.44 & \phantom{11}140.01\\
\cline{2-8}
& {$7$} &$t_f$  & \phantom{1}17.74 (3)  &	 \phantom{1}26.56 (3)  &	 \phantom{11}56.25 (3) & \phantom{11}281.27 (3) & 6.14\\
& {$8$} &$t_f$  & \phantom{1}48.39 (2)  &	 \phantom{1}74.40 (2)  &	\phantom{1}153.35 (2) & \phantom{11}506.84 (2) & 33.25\\
& {$9$} &$t_f$  & 281.27 (3) &	462.52 (3)  &1184.74 (3) & \phantom{1}6261.34 (3) & 342.98\\
\hline
\end{tabular}
\end{center}
\end{table}

When the viscosity parameter is small (i.e., $\nu=1/600$), the coarse-grid computation requires the $\kappa$-term approximation on a relatively fine spatial grid (i.e., $\ell = 6$).  The exponential boundary layer gets narrower as the viscosity parameter gets smaller, which requires the use of a finer spatial grid for the coarse-grid computation. If the coarse-grid computation is performed on coarser spatial grids, it fails to identify the rank structure of solutions and to yield a proper truncation operator. Analogously, when the number of terms, $M$, in the KL expansion \eqref{kl_exp} is large, the coarse-grid computation has to be done on a relatively fine spatial grid because the KL expansion contains more spatially oscillatory terms. In the last columns of Table \ref{tab_cd1} and \ref{tab_cd2}, $\ast$ and $\dagger$ indicate that the coarse-grid solutions are computed on the level 5 and the level 6 spatial grid, respectively. 

\vspace{0.5mm}
\textbf{Comparison to a truncation operator based on singular values.} We again compare the performance of the proposed solver to the preconditioned low-rank projection method combined with the conventional truncation operator, the SVD-based truncation operator. Table \ref{tab-svd-cd} shows the computation time required to compute approximate solutions using the conventional and the new truncation strategy.  When the low-rank projection method with SVD-based truncation operator is used, initial values for rank $\kappa$ in Algorithm \ref{lrp_alg} are obtained from coarse-grid computations of the multilevel rank reduction strategy. 

\begin{table}[htbp]
\caption{CPU times to compute approximate solutions  of the convection-diffusion equation satisfying $\epsilon = 10^{-5}$ and $10^{-6}$ using the preconditioned low-rank projection (LRP) methods with the multilevel rank-reduction and the singular value based truncation on the level 8 spatial grid (i.e., $n_x = 257^2$).}
\begin{center}\footnotesize
\renewcommand{\arraystretch}{1.3}
\label{tab-svd-cd}
\begin{tabular}{c| c |  c | c |rrrr }
 \hline
&Viscosity ($\nu$) &Solver & &  $M=5$  & $M=7$  & $M=10$ & $M=15$  \\
 \hline
\multirow{10}{*}{$\epsilon = 10^{-5}$} & \multirow{2}{*}{$1/20$}&LRP-SVD & $t_{\text{SVD}}$& 68.45 & 100.83 & 201.34 & 438.25 \\
\cline{3-8}
&&LRP-Multilevel &$t$& 54.06  & 72.08 & 154.79 & 338.21 \\
\cline{2-8}
&\multirow{2}{*}{$1/100$}&LRP-SVD & $t_{\text{SVD}}$& 93.91 & 121.89 & 295.27 & 655.71 \\
\cline{3-8}
&&LRP-Multilevel &$t$& 55.28 & 64.36 & 125.88 & 285.94 \\
\cline{2-8}
&\multirow{2}{*}{$1/200$}&LRP-SVD & $t_{\text{SVD}}$& 90.70 & 122.56 & 251.60 & 574.68 \\
\cline{3-8}
&&LRP-Multilevel &$t$& 55.42 & 66.08 & 115.68 & 258.97 \\
\cline{2-8}
&\multirow{2}{*}{$1/400$}&LRP-SVD & $t_{\text{SVD}}$& 91.11 & 107.47 & 221.32 & 475.60 \\
\cline{3-8}
&&LRP-Multilevel &$t$& 69.01  & 76.63 & 158.07 & 416.36 \\
\cline{2-8}
&\multirow{2}{*}{$ 1/600$}&LRP-SVD & $t_{\text{SVD}}$& 90.33  & 103.44 & 218.35 & 484.08\\
\cline{3-8}
&&LRP-Multilevel &$t$& 75.26  & 86.48 & 176.93 & 422.85\\
\hline
\multirow{10}{*}{$\epsilon = 10^{-6}$}&\multirow{2}{*}{$1/20$}&LRP-SVD & $t_{\text{SVD}}$& 132.08 & 234.15 & 570.56 & 1748.43\\
\cline{3-8}
&&LRP-Multilevel &$t$& 86.74 & 145.86 & 401.83 & 1267.71 \\
\cline{2-8}
&\multirow{2}{*}{$1/100$}&LRP-SVD & $t_{\text{SVD}}$& 121.88 & 196.66 & 471.11 & 1479.80\\
\cline{3-8}
&&LRP-Multilevel &$t$& 84.97 & 126.77 & 329.52 & 1088.05\\
\cline{2-8}
&\multirow{2}{*}{$1/200$}&LRP-SVD & $t_{\text{SVD}}$& 106.79 & 188.76 & 416.52 & 1203.78\\
\cline{3-8}
&&LRP-Multilevel &$t$& 77.79 & 128.96 & 293.30 & 892.18\\
\cline{2-8}
&\multirow{2}{*}{$1/400$}&LRP-SVD & $t_{\text{SVD}}$& 107.12 & 168.01 & 380.01 & 1015.88\\
\cline{3-8}
&&LRP-Multilevel &$t$& 78.04 & 112.55 & 269.48 & 797.50 \\
\cline{2-8}
&\multirow{2}{*}{$ 1/600$}&LRP-SVD & $t_{\text{SVD}}$& 122.44  & 231.07 & 421.76 & 1208.88 \\
\cline{3-8}
&&LRP-Multilevel &$t$& 97.00  & 129.87 & 234.00 & 670.90\\
\hline
\end{tabular}
\end{center}
\end{table}

\subsection{Choices of coarse spatial grid} \label{sec_coarse_sel}
Finally, we discuss criteria for choosing the coarse grid used to generate truncation operators. The basic idea is that the coarse grid needs to be fine enough so that important features of the problem are represented. This quality is problem dependent, and we outline what is needed for the two 
types of problems we examined.

First consider the diffusion equation of Section 5.1. Here, the issue is the oscillatory nature of components of the random field $a(\pmb{x},\,\xi)$. In the the KL expansion \eqref{kl_exp}, the eigenpairs, $\{(\lambda_i,\,a_i(\pmb{x}))\}_{i=1}^M$, can be obtained by solving the following integral equation, 
\begin{equation}
\label{eig_int}
\int_D C(\pmb{x},\,\pmb{y}) a_i(\pmb{y}) d\pmb{y} = \lambda_i a_i(\pmb{x}), \qquad i=1,\,\dots,\,M
\end{equation}
where $C(\pmb{x},\,\pmb{y})$ is the covariance kernel \eqref{exp_cov}. Since the kernel is separable, the eigenfunctions of the integral problem \eqref{eig_int} can be represented as $a_i(\pmb{x}) = a_k^1(x_1) a_j^2 (x_2)$, where $\{a_k^1\}_{k=1}^{\infty}$ and $\{a_j^2\}_{j=1}^{\infty}$ are the eigenfunctions of the one-dimensional integral problem (i.e., $\int_D \text{exp}(-| x_l - y_l |/c) a_k^l(y_l) dy_l = \lambda_k^l a_k^l(x_l)$, $l=1,\,2$). The eigenvalues, $\{ \lambda_i \}_{i=1}^M$, are in decreasing order and $\lambda_i$ is the $i$th largest value of products $\lambda_k^1 \lambda_j^2$ for $k,\,j=1,\,2,\,\cdots$. Analytic expressions for the 1D eigenfunctions are given in \cite{ghanem2003stochastic} as, for $l=1,\,2$, 
\begin{align}
\begin{split}
\label{eig_fun}
a_k^l(x) &= h^l(\theta_k) \cos (\theta_k x) \qquad\: \text{for even } k, \\
{a_k^{l\ast}}(x) &= h^{l\ast}(\theta_k^\ast) \sin (\theta_k^\ast x) \qquad \text{for odd } k,
\end{split}
\end{align}
where $\theta_k$ and $\theta^\ast_k$ are the solutions of 
\begin{equation*}
\frac{1}{c} - \theta \tan \left( \frac{\theta}{2} \right) = 0 \quad \text{ and } \quad \theta^\ast + \frac{1}{c} \tan \left(\frac{\theta^\ast}{2}\right) =0,
\end{equation*}
respectively, when the 1D integral problem is posed on $[-\frac{1}{2},\,\frac{1}{2}]$. As $i$ in the KL expansion \eqref{kl_exp} increases, the eigenfunctions $a_i(\pmb{x})$ become more oscillatory over the spatial domain (i.e., $\theta_k$ or $\theta_k^\ast$ become larger), so that finer coarse spatial grids are required to capture the oscillatory features of the KL expansion. Table \ref{tab_freq} shows the largest value of $\{\theta_k,\, \theta_k^\ast\}$ of the eigenfunctions in the KL expansion, the half-wavelength of 
the functions from \eqref{eig_fun} and our choice of coarse spatial grid refinement levels, $\ell^c$, for different values of $M$. With these coarse grids, there are approximately eight grid points per half wave, enough to capture the qualitative character of the wave.

\begin{table}[tbp]
\caption{Largest values of $\theta_k$ or $\theta_k^\ast$ of eigenfunctions \eqref{eig_fun} in the KL expansion, required grid refinement level $\ell^c$, half wavelength $\pi/\theta$, and element size $h^c = 2^{-\ell^c}$ for different values of $M$.}
\begin{center}\footnotesize
\renewcommand{\arraystretch}{1.3}
\label{tab_freq}
\begin{tabular}{ c |  cccc cc }
\hline
$M$ & 3 & 5 & 7 & 10 & 15 & 20\\
\hline
$\underset{}{\operatorname{max}}(\theta_k,\,\theta_k^\ast)$ & 3.25 & 6.36 & 9.49 & 12.63 & 18.90 & 25.19\\
wavelength/2 & \phantom{1}.97 & \phantom{1}.49 & \phantom{1}.33 & \phantom{11}.25 & \phantom{11}.17 & \phantom{11}.12\\
$\ell^c$ $(h^c)$  & 3 $\left(\frac{1}{8}\right)$& 4 $\left(\frac{1}{16}\right)$& 4 $\left(\frac{1}{16}\right)$ & 5 $\left(\frac{1}{32}\right)$& 5 $\left(\frac{1}{32}\right)$& 6 $\left(\frac{1}{64}\right)$\\
\hline
\end{tabular}
\end{center}
\end{table}

\begin{figure}[!tbp]
	\centering
	\includegraphics[angle=0, scale=0.45]{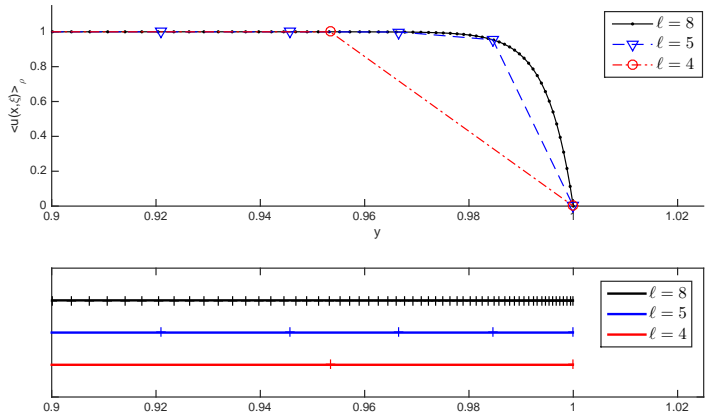}
	\caption{Mean solutions $\langle u(\pmb{x},\,\xi) \rangle_\rho$ at $x = 1$ and $y=[0.9,\,1]$ illustrating the exponential boundary layer for varying spatial grid refinement level, $\ell = \{4,\,5,\,8\}$, (top) and lengths in $y$-direction of first few elements from $y=1$ (bottom).}
	\label{exp_layer}
\end{figure}

We turn now to the convection-diffusion equation of Section 5.2. This problem has the same diffusion coefficient \eqref{kl_exp} as the diffusion problem, but in addition its solution has an exponential boundary layer. In particular, for small $\nu$, the width of the layer is smaller than the finest interval needed to represent the eigenfunctions in \eqref{kl_exp}, and in this case the coarse grid must be finer than that needed for the diffusion problem (whose solution is smooth). In Figure \ref{exp_layer}, the top plot illustrates the mean solutions $\langle u(\pmb{x},\,\xi) \rangle_\rho$ of the weak formulation of \eqref{cd_strong} at $x=1$, which are computed on two coarse spatial grids $\ell = \{4,\,5\}$ using PGD and a fine spatial grid $\ell = 8$ using the multilevel method, with the viscosity parameter, $\nu = \frac{1}{200}$, and $M=10$ random variables. The bottom plot shows the lengths of the first few elements in the $y$-direction near $y=1$ for these refinement levels. If the level-4 spatial grid is used for the coarse grid computation (i.e., $\ell=4$, red line in Figure \ref{exp_layer}), the width of exponential boundary layer is much narrower than the length of the smallest element and the coarse-grid solution gives a poor representation of the boundary layer. When this coarse grid is used to construct the truncation operator, the multilevel scheme fails to compute an accurate approximate solution on a fine spatial grid (i.e., $\ell=8$, black line in Figure \ref{exp_layer}). On the other hand, the level-5 spatial grid (i.e., $\ell=5$, blue line in Figure \ref{exp_layer}) is fine enough for the coarse-grid solution to represent the character of the exponential layer, and with this coarse-grid, the resulting multilevel scheme efficiently computes an accurate fine-grid solution.

Although this discussion shows that some a priori knowledge of the problem is needed to identify the coarse grid operator, in general this information is not difficult to come by. In particular, we are assuming that the expansion (2.5) is known, and it is straightforward to identify the resolution needed to represent its components, for example by examining one-dimensional cross-sections of them. If as for the second problem some knowledge of the solution is needed, this can be obtained cheaply from the solution of a deterministic problem derived from the mean of the diffusion coefficient; indeed, for the convection-diffusion problem, the boundary layer for the deterministic solution has essentially the same character as that of the stochastic solution whose mean is shown in Figure 2.

\section{Conclusion} \label{sec_con}
We have studied iterative solvers for low-rank solutions of sto-chastic Galerkin systems of stochastic partial differential equations. In particular, we have explored low-rank projection methods in tensor format for linear systems of Kronecker-product structure. For the computational efficiency of the projection methods, basis vectors and iterates in the projection methods are forced to have low rank, which is achieved by a multilevel rank-reduction strategy. We have examined the performance of this strategy with two benchmark problems: stochastic diffusion problems and stochastic convection-diffusion problems. For both problem classes, the rank structure of the solution can be identified by an inexpensive coarse-grid computation, and with the resulting multilevel rank-reduction strategy, the low-rank projection method is more efficient than methods for which the truncation operator is based on singular values.

\bibliography{mybib}
\bibliographystyle{siam}

\end{document}